 \documentstyle[graphicx,amscd,amssymb,verbatim,12pt,righttag,color]{amsart}
 \newtheorem{theorem}{Theorem}[section]
 
 \newtheorem{Prop}[theorem]{Proposition}
 \newtheorem{Lem}[theorem]{Lemma}
 
 \newtheorem{Rem}[theorem]{Remark}
 \newtheorem{Example}[theorem]{Example}

\newcommand{\eproof}{\hfill$\square$}%\rule{2.2mm}{3.0mm}}

 \numberwithin{equation}{section}
 \renewcommand{\rm}{\normalshape}
\begin{document}

\title  {Frames of multi-windowed exponentials on subsets of ${\mathbb R}^d$}

%\date{29 November 2011}
\author{Jean-Pierre Gabardo}
\address{Department of Mathematics and Statistics, McMaster University,
Hamilton, Ontario, L8S 4K1, Canada}
\email{gabardo@@mcmaster.ca}
\author{Chun-Kit Lai}
\address{Department of Mathematics and Statistics, McMaster University,
Hamilton, Ontario, L8S 4K1, Canada}
\email{cklai@@math.mcmaster.ca}

%\thanks{The research is partially supported by the RGC grant of Hong Kong and the Focused
%Investment Scheme of CUHK; the first two authors are also supported by the National Natural Science Foundation of
%China 10871180.}
%
\date{\today}
 %Secondary 42C15.}
%%\thanks{ The research is supported in part by the HKRGC Grant and the Focused Investment Scheme of CUHK. }
\maketitle

\begin{abstract}
Given discrete subsets $\Lambda_j\subset {\Bbb R}^d$, $j=1,\dots,q$,
consider the set of windowed exponentials
$\bigcup_{j=1}^{q}\{g_j(x)e^{2\pi i \langle\lambda,x\rangle}: \lambda\in\Lambda_j\}$ on $L^2(\Omega)$.
We show that a necessary and sufficient condition for the windows $g_j$ to form a
frame of  windowed exponentials for $L^2(\Omega)$ with some $\Lambda_j$ is
that $m\leq \max_{j\in J}|g_j|\leq M$ almost everywhere on $\Omega$ for some
subset $J$ of $\{1,\cdots, q\}$. If $\Omega$ is unbounded, we show that there is
no frame of windowed exponentials if the Lebesgue measure of $\Omega$ is infinite.
If $\Omega$ is unbounded but of finite measure, we give a sufficient condition for
the existence of Fourier frames on $L^2(\Omega)$. At the same time, we also construct
examples of unbounded sets with finite measure that have no tight exponential frame.
\end{abstract}

\section{introduction}

Let $\Omega$ be a Lebesgue measurable set on ${\Bbb R}^d$ and let $g_j\in
 L^2(\Omega)\setminus\{0\}$, $j=1,\cdots, q$ and $q<\infty$.
 Let also $\Lambda_j$, $j=1,\cdots, q$ be some countable sets on ${\Bbb R}^d$.
The collection $\bigcup_{j=1}^{q}{\mathcal E}(g_j,\Lambda_j) =
\bigcup_{j=1}^{q}\{g_j(x)e^{2\pi i \langle\lambda,x\rangle}: \lambda\in\Lambda_j \}$
is called a set of {\it windowed exponentials} with {\it windows} $g_j$.
Recall that $\bigcup_{j=1}^{q}{\mathcal E}(g_j,\Lambda_j)$ is a {\it frame} for $L^2(\Omega)$ if
there exist $A,B>0$ such that
\begin{equation}\label{eq0.1}
A\|f\|_{L^2(\Omega)}^2\leq \sum_{j=1}^{q}\sum_{\lambda\in\Lambda_j}
\left|\int_{\Omega}f(x)\overline{g_j(x)}e^{-2\pi i \langle\lambda,x\rangle}dx\right|^2
\leq B\|f\|_{L^2(\Omega)}^2
\end{equation}
for all $f\in L^2(\Omega)$, where $\|f\|_{L^2(\Omega)}^2=\int_{\Omega}|f(x)|^2dx$.
If the second inequality in (\ref{eq0.1}) is satisfied, then the set of the windowed
exponentials is called a {\it Bessel sequence.} If the collection is generated by the
single window  $g= \chi_{\Omega}$, it is called a  {\it Fourier frame}.

\medskip

 The study of Fourier frames was initiated by Duffin and Schaeffer in their work on non-harmonic
 Fourier series \cite{[DS]}. The existence of Fourier frames
$\{e^{2\pi i \langle\lambda,\cdot\rangle}\}_{\lambda\in\Lambda}$ on $L^2(\Omega)$
 was also known to be equivalent to the sampling problems on the
Paley-Wiener space $PW_{\Omega}$, which ask for the reconstruction of
the band-limited functions $f$ by their sampled values $\{f(\lambda)\}$ (see \cite{[Y]}).
Nowadays, Fourier frames and, more generally,
windowed exponentials  have a wide range of applications in
different area of mathematics, engineering and signal processing \cite{[AG],[Chr]}.

\medskip

Windowed exponentials also arise naturally in frame theory.
 Gr\"{o}chenig and Razafinjatovo \cite{[GR]} derived the
famous necessary Beurling density condition of Landau \cite{[Lan]}
on Fourier frames by considering windowed exponentials.
It is also known that the study of frame of translates and regular Gabor frames can be
reduced to that of windowed exponentials via the Fourier transform and the
Zak transform respectively \cite{[Chr],[G]}. Heil et al have
recently made extensive studies in the basis properties,
density conditions and different aspects of the windowed exponentials \cite{[HK],[HY]},
and the reader can refer to \cite{[H]} for a comprehensive introduction to the theory
 of  windowed exponentials.

\medskip

 In this paper, we will give a complete characterization of the collections of windows $g_j$
with the property that $\bigcup_{j=1}^{q}{\mathcal E}(g_j,\Lambda_j)$
form a frame on $L^2(\Omega)$ for some discrete sets $\Lambda_j\subset {\Bbb R}^d$.
This characterization is motivated by the recent work of the second named author on
Fourier frames of absolutely continuous measures \cite{[Lai],[DL]}.
In fact, we will see that Fourier frames of absolutely continuous measures
are equivalent to frames of windowed exponentials generated by a single window
(Proposition \ref{prop4.1}). Therefore, our result will be a further generalization.
 Denoting by $|\Omega|$ the Lebesgue measure of $\Omega$, we have the following result.

 \begin{theorem}\label{th0.1}
 If $|\Omega|$ is infinite, then  there is no frame of windowed exponentials on $L^2(\Omega)$.
\end{theorem}

This statement is no longer true if we allow infinitely many windows.
For example on ${\Bbb R}^d$, the system
$$
\bigcup_{n\in{\Bbb Z}^d}\{\chi_{[0,1]^d+n}(\cdot)
e^{2\pi i \langle m, \cdot\rangle}: m\in{\Bbb Z}^d\}  =
\{e^{2\pi i \langle m, \cdot\rangle}\chi_{[0,1]^d}(\cdot-n):m,n\in{\Bbb Z}^d\}.
$$
is a  standard example of a Gabor orthonormal basis for $L^2({\Bbb R}^d)$.
\medskip

If the measure of $\Omega$ is finite, we need to separate our
analysis into the case where $\Omega$ is bounded or unbounded.

\begin{theorem}\label{th0.2}
Let $\Omega\subset{\mathbb R}^d$ be a bounded Lebesgue measurable set and let
$g_j$, $j=1,2\cdots, q$, be a finite set of functions in $L^2(\Omega)$.  Let also
$$J = \{j: \|g_j\|_{\infty}<\infty\}.
$$
Then there exists $\Lambda_j$ such that  $\bigcup_{j=1}^{q}{\mathcal E}(g_j,\Lambda_j)$
form a frame in $L^2(\Omega)$ if and only if there exists  $m>0$ such that
$$
\max_{j\in J}|g_j|\geq m
$$
almost everywhere on $\Omega$.
\end{theorem}

An intermediate result of independent interest, Theorem \ref{th3.1}, is needed in
the proof of the above theorem, in which explicit upper and lower bounds for
the quantities $\max_{j\in J}|g_j|$ are given  in terms of the frame bounds and
upper Beurling densities of certain measures associated with the sets $\Lambda_j$.

\medskip

The characterization given in Theorem \ref{th0.2}  implies that the unbounded functions in the
original collection are actually not needed in producing frame of windowed exponentials
and we just need to check whether the maximum of the moduli of the
remaining bounded functions is bounded away from $0$ a.e.  on $\Omega$. In particular,
if the collection consists only of unbounded functions, it cannot form any frame of windowed exponentials.

\medskip

One essential ingredient in our proofs is a surprising relationship between
Beurling densities and the bounds in some convolution inequalities developed in
\cite{[Ga1],[Ga2]}. Convolution inequalities arise naturally in the study of frame theory,
tilings and spectral sets. Making use of this relationship has the advantage to simplify many
technical calculations and to allow the theorems above
to hold in the more general setting of {\it generalized frames of windowed exponentials}
(Remark \ref{rem4.1}) without much additional work.

\medskip

The situation unfortunately becomes vastly more complicated
if the set $\Omega$ is unbounded but  still of finite Lebesgue measure.
The necessary condition given in Theorem \ref{th0.2} for a system of windowed exponentials
to form a frame for $L^2(\Omega)$ still holds in the unbounded case,
as the proof only uses the fact that $\Omega$ has finite Lebesgue measure
(see Theorem \ref{th3.1}). However, we do not know of any
unified argument to show that frames of windowed exponentials for $L^2(\Omega)$ exist for every
such $\Omega$. In all the known examples,
the construction of frames of windowed exponentials and Fourier frames
is based on a tight frame defined on a larger set.
We therefore examine the existence of tight Fourier frames for unbounded sets of finite measure.
It is not hard to  prove that {\it if there exists a lattice $\Gamma$ such that
$\sum_{\gamma\in\Gamma}\chi_{\Omega}(\cdot+\gamma)\leq 1$
(i.e. elements in $\Omega$ are distinct residue class of $\Gamma$),
then $L^2(\Omega)$ will admit a tight Fourier frame.} (Proposition \ref{prop5.1}).
We don't know whether or not this condition is necessary  but,
on the other extreme, we can construct examples where no
 tight frames can exist using the following theorem.

\begin{theorem}\label{th0.3}
Suppose $\Omega$ is a measurable set of finite Lebesgue measure such that
$|\Omega\cap \Omega+x|>0$ for all $x\in{\Bbb R}^d$ with $|x|>R $ for some $R>0$,
then $L^2(\Omega)$ does not admit any tight Fourier frame
(i.e. a Fourier frame with $A=B$ in (\ref{eq0.1})).
\end{theorem}

 Examples of sets satisfying the conditions in the above theorem are
 not difficult to obtain and the theorem shows that ordinary method of
 Fourier frame construction fails for these. However, we cannot prove
 whether or not Fourier frames always exist for such sets.

\medskip

If the condition $|\Omega\cap (\Omega+x)|>0$ is valid for all
 $x\in {\Bbb R}^d$, we can strengthen the conclusion of Theorem \ref{th0.3} to
obtain that the only {\it tight frame measures} (see (\ref{eq4.2}))
 for $L^2(\Omega)$ are the positive multiples of the Lebesgue measure
on ${\Bbb R}^d$ (Theorem \ref{th5.1}).

\medskip

We organize the paper as follows. We will give some preliminaries
 on convolution inequalities in Section 2. We then prove Theorem \ref{th0.1}
and \ref{th0.2} in Section 3. After that, we discuss frames on unbounded sets
of finite measure and prove Theorem \ref{th0.3} in Section 4.
In the last section, we will apply our results on windowed exponentials
to related systems: frames of translates, frames of absolutely continuous measures
and Gabor frames. Although these results are known, we are able to
recover them with a new approach and simpler proofs.

\bigskip

\section{Preliminaries}

Let $\mu$ be a positive Borel measure on ${\Bbb R}^d$. We define its associated
{\it upper and lower Beurling densities} as
$$
D^{+}(\mu) = \limsup_{h\rightarrow\infty}\sup_{x\in{\Bbb R}^d}\frac{\mu(x+Q_{h})}{h^d}, \  D^{-}(\mu)
= \liminf_{h\rightarrow\infty}\inf_{x\in{\Bbb R}^d}\frac{\mu(x+Q_{h})}{h^d},
$$
where $Q_h$ is the cube of side length $h$ centered at the origin.
If $\Lambda$ is a countable set on ${\Bbb R}^d$ and we denote by
$\delta_{\Lambda}$ the measure $\sum_{\lambda\in\Lambda}\delta_{\lambda}$,
then the above definitions of  Beurling densities become those of the usual Beurling
densities of a  discrete set $\Lambda$, which we denote by $D^{+}(\Lambda)$ and
$D^{-}(\Lambda)$, respectively.
We say that a measure $\mu$ is {\it translation-bounded} if for every
compact set $K$ there exists a constant $C_K>0$ such that $\mu(x+K)\leq C_K$ for all $x\in{\Bbb R}^d$.

\medskip

From the definition above, we can easily obtain the following inequalities.

\begin{Prop}\label{prop2.1}
Let $\mu$ and $\nu$ be  positive Borel measures on ${\Bbb R}^d$. Then
$$
D^{+}(\mu)\leq D^{+}(\mu+\nu)\leq D^{+}(\mu)+D^{+}(\nu).
$$
In particular, if  $D^{+}(\nu)=0$, then $D^{+}(\mu+\nu)= D^{+}(\mu)$.
\end{Prop}

\begin{pf}
 We note that from the definition, we immediately have
$$
\sup_{x\in{\Bbb R}^d}\frac{\mu(x+Q_h)}{h^d}\leq\sup_{x\in{\Bbb R}^d}
\frac{(\mu+\nu)(x+Q_h)}{h^d}\leq \sup_{x\in{\Bbb R}^d}
\frac{\mu(x+Q_h)}{h^d}+\sup_{x\in{\Bbb R}^d}\frac{\nu(x+Q_h)}{h^d}.
$$
Hence, passing to the limit, we have $D^{+}(\mu)\leq
 D^{+}({\mu+\nu})\leq D^{+}({\mu})+D^{+}(\nu)$. The second statement is clear from the inequalities.
\end{pf}

\begin{Rem} {\rm It should be pointed out that Proposition \ref{prop2.1} is not
true for the lower Beurling density. To see this, we can let
$\mu = \sum_{n=0}^{\infty}\delta_n$ and $\nu=\sum_{n=1}^{\infty}\delta_{-n}$.
Then $D^{-}(\mu) = D^{-}(\nu)=0$, but $D^{-}(\mu+\nu) = D^{-}({\Bbb Z}) =1$.}
\end{Rem}

\medskip

  For a positive Borel measure $\mu$ and an locally integrable function $f\ge 0$
on ${\Bbb R}^d$,
we define the convolution $\mu\ast f$ using the formula
   $$
\mu\ast f(x): = \int f(x-y)d\mu(y), \ x\in{\mathbb R}^d.
$$
In our proofs, we need to exploit a relationship between the convolutions and
the Beurling densities for several measures. The following theorem was
 proved in \cite[Corollary 6 and 7]{[Ga1]}.
%It is also convenient to define the set of all probability densities.
%$$
%{\mathcal P}({\Bbb R}^d) = \{f: f\geq 0, \ \int_{\Bbb R^d} f(x)dx=1\}
%$$

%\begin{theorem}\label{th2.1}
%Let $\mu$ be a positive Borel measure on ${\Bbb R}^d$. Then
%
%\vspace{0.2cm}
%
%\noindent{\rm (i)} $D^{+}(\mu) = \inf\{C: \mu\ast f\leq C \ \mbox{ {\rm a.e. for some} } \ f\in{\mathcal P}({\Bbb R}^d) \}$.
%
%\vspace{0.2cm}
%
%\noindent{\rm (ii)} If $\mu$ is translation-bounded, then $D^{-}(\mu) = \sup\{D: \mu\ast f\geq D \ \mbox{ {\rm a.e. for some} } \ f\in{\mathcal P}({\Bbb R}^d) \}$.
%\end{theorem}
%
%An immediate corollary of Theorem \ref{th2.1} is
%
%\begin{Cor}\label{th2.2}
%Let $\mu$ be a positive translation-bounded Borel measure on ${\Bbb R}^d$.
%
%\vspace{0.2cm}
%
%\noindent{\rm (i).} Suppose that there exists $h\in L^1({\Bbb R}^d)$ and $B>0$ such that $\mu\ast h\leq B$ almost everywhere.  Then  $D^{+}(\mu)\int_{{\Bbb R}^d}h(x)dx\leq B.$
%
%\vspace{0.2cm}
%
%\noindent{\rm (ii).} Suppose that there exists $h\in L^1({\Bbb R}^d)$ and $A>0$ such that $A\leq\mu\ast h\leq B$ almost everywhere.  Then  $A\leq D^{-}(\mu)\int_{{\Bbb R}^d}h(x)dx$.
%\end{Cor}

\begin{theorem}\label{th2.3}
For $i=1,\cdots,m$, let $\mu_i$  be positive Borel measures on ${\Bbb R}^d$,
let $h_i$ be non-negative functions in $L^1({\Bbb R}^d)$
and write  $\mu =\sum_{i=1}^{m}(\int_{{\Bbb R}^d} h_i(x)\,dx)\,\mu_i$.

\vspace{0.2cm}

\noindent{\rm (i)} Suppose that there exists $B>0$ such that
$\sum_{i=1}^{m}\mu_i\ast h_i\leq B \ \mbox{a.e. on } \ {\Bbb R}^d$, then
$D^{+}(\mu)\leq B$.

\vspace{0.2cm}

\noindent{\rm(ii)} Suppose that there exists $A>0$ such that
$A\leq \sum_{i=1}^{m}\mu_i\ast h_i \ \mbox{a.e. on } \ {\Bbb R}^d$ and
that all the $\mu_i$ are translation-bounded, then
$A\leq  D^{-}(\mu)$.
\end{theorem}

We also recall an important condition equivalent to the translation-boundedness
of a measure $\mu$ \cite[Proposition 1]{[Ga1]}.

\begin{Prop}\label{Prop2.2-}
Let $\mu$ be a positive Borel measure measure on ${\Bbb R}^d$.
Then $\mu$ is translation-bounded if and only if there exists $f\in L^1({\Bbb R}^d)$  with $f\ge 0$
and a constant $C>0$ such that $\mu\ast f\leq C$ a.e. on ${\Bbb R}^d$.
\end{Prop}

The Fourier transform a function $f\in L^1({\Bbb R}^d)$ is defined by
 $$
 \widehat{f}(\xi) = \int_{{\Bbb R}^d} f(x)e^{-2\pi i \langle\xi,x\rangle}dx,\quad \xi\in {\Bbb R}^d,
 $$
and extended in the usual way as a unitary operator on $f\in L^2({\Bbb R}^d)$.
  The following proposition illustrates how  convolution inequalities appear
naturally when dealing with Bessel systems or frames of windowed exponentials.
 \medskip

\begin{Prop}\label{prop2.2}

\vspace{0.2cm}

\noindent{(i)}\
Let $\bigcup_{j=1}^{q}{\mathcal E}(g_j,\Lambda_j)$ be a
Bessel sequence of windowed exponentials on $L^2(\Omega)$
(where $|\Omega|$ can be finite or infinite). Then, for any $f\in L^2(\Omega)$
such that $fg_j\in L^2(\Omega)$ for all $j$, we have
 $$
 D^{+}(\mu_f)\leq B\|f\|_{L^2(\Omega)}^2
 $$
where $\mu_f = \sum_{j=1}^{q}\left(\int_{\Omega}|fg_j|^2\right)\,\delta_{\Lambda_j}$
and  all the measures $\delta_{\Lambda_j}$, $j=1,\dots, q$, are translation-bounded.
\vspace{0.2cm}

\noindent{(ii)} \
If, furthermore, the collection  $\bigcup_{j=1}^{q}{\mathcal E}(g_j,\mu_j)$ is
a frame of windowed exponentials for $L^2(\Omega)$,
then $$A\|f\|_{L^2(\Omega)}^2\leq D^{-}(\mu_f).$$
\end{Prop}

\begin{pf}
(i) Replacing $f$ with the function  $\overline{f}e^{2\pi i \langle\xi,\cdot\rangle}$
in the definition of Bessel sequence of windowed exponentials in (\ref{eq0.1}), we have
$$
 \sum_{j=1}^{q}\sum_{\lambda\in\Lambda_j}\left|\widehat{\chi_{\Omega}fg_j}(\xi-\lambda)\right|^2
\leq B\|f\|_{L^2(\Omega)}^2.
$$
In particular, we can write the inner sum as
$$
\sum_{\lambda\in\Lambda_j}\left|\widehat{\chi_{\Omega}fg_j}(\xi-\lambda)\right|^2 =
 \delta_{\Lambda_j}\ast (|\widehat{\chi_{\Omega}fg_j}|^2)(\xi)
$$
to obtain
\begin{equation}\label{eq2.1}
 \sum_{j=1}^{q}\delta_{\Lambda_j}\ast (|\widehat{\chi_{\Omega}fg_j}|^2)(\xi)\leq B\|f\|_{L^2(\Omega)}^2.
\end{equation}
As $\int|\widehat{\chi_{\Omega}fg_j}|^2=\int_{\Omega}|fg_j|^2(<\infty)$
from the Plancherel identity and the assumption on $f$, it follows from Theorem \ref{th2.3} (i) that
$$
  D^{+}(\mu_f) \leq B\|f\|_{L^2(\Omega)}^2,
$$
where $\mu_f = \sum_{j=1}^{q}\left(\int_{\Omega}|fg_j|^2\right)\,\delta_{\Lambda_j}$.
From (\ref{eq2.1}), for all $j=1,\cdots,q$, we have
$$
\delta_{\Lambda_j}\ast(|\widehat{\chi_{\Omega}fg_j}|^2)(\xi)
\leq B\|f\|_{L^2(\Omega)}^2,\quad \xi\in{\mathbb R}^d.
$$
and thus each measure $\delta_{\Lambda_j}$ is translation-bounded by Proposition \ref{Prop2.2-}.

\medskip

(ii) By (i), all the measures $\delta_{\Lambda_j}$ are translation-bounded.
By an argument similar to the one used in (i), we obtain
$$
A\|f\|_{L^2(\Omega)}^2\leq\sum_{j=1}^{q}\delta_{\Lambda_j}\ast(|\widehat{\chi_{\Omega}fg_j}|^2)(\xi).
$$
By Theorem \ref{th2.3} (ii), the conclusion follows.
\end{pf}

\bigskip

\section{windowed exponentials}

In this section, we will prove our main results concerning
general frames of windowed exponentials.
We first need a lemma.

\begin{Lem}\label{lem3.1}
Let $\bigcup_{j=1}^{q}{\mathcal E}(g_j,\Lambda_j)$ be a Bessel sequence of windowed
exponentials for $L^2(\Omega)$ with $g_j\neq 0$ for all $j$.
\vspace{0.2cm}

\noindent{(i)} If $|\Omega|<\infty$, then $D^{+}\Lambda_j<\infty$ for all $j$ and
at least one of the $\Lambda_j$ has positive upper Beurling density if the windowed
exponentials form a frame for $L^2(\Omega)$.

\vspace{0.2cm}

\noindent{(ii)} If the windowed exponentials form a frame for $L^2(\Omega)$ and
$|\Omega|=\infty$, then  we have $D^{+}\Lambda_j=\infty$ for all $j$.

\end{Lem}

\begin{pf}
(i) If $|\Omega|$ is finite, letting $f=\chi_{\Omega}$ in Proposition \ref{prop2.2},
we have
$$
D^{+}(\mu_f)\leq B|\Omega|,
$$
where $\mu_f = \sum_{j=1}^{q}\|g_j\|_{L^2(\Omega)}^2\delta_{\Lambda_j}$.
Letting $m = \min_{j}\,\|g_j\|_{L^2(\Omega)}^2>0$, we have the inequality $\mu_f\geq m\sum_{j=1}^{q}
\delta_{\Lambda_j}$. Hence, invoking Proposition \ref{prop2.1}, it follows that
$$
D^{+}(\Lambda_j)\leq D^{+}(\sum_{i=1}^{q}\delta_{\Lambda_i})\leq\frac{1}{m}D^{+}(\mu_f)
\leq \frac{B|\Omega|}{m}<\infty,\quad j=1,\dots, q.
$$
In addition, if the set of the windowed exponentials is a frame, we have also the inequality
 $D^{-}(\mu_f)\geq A|\Omega|.$ Letting $M = \max_{j}\,\|g_j\|_{L^2(\Omega)}^2$, we have
\begin{equation}\label{eq3.1}
0<\frac{A|\Omega|}{M}\leq \frac{1}{M}D^{+}(\mu_f)\leq D^{+}(\sum_{j=1}^{q}\delta_{\Lambda_j})
\leq \sum_{j=1}^{q}D^{+}(\Lambda_j),
\end{equation}
showing that $D^{+}(\Lambda_j)>0$ for some $j$.

\medskip

\noindent{(ii).} If $|\Omega|=\infty$, let $\Omega_N=\Omega\cap Q_N$,
where $Q_N$ is the cube of side length $N$ centered at origin.
Then $L^2(\Omega_N)\subset L^2(\Omega)$ and it is easy to see that
$\bigcup_{j=1}^{q}{\mathcal E}(g_j,\Lambda_j)$ is still a frame of $L^2(\Omega_N)$
for all $N$ large enough so that $|\Omega_N|>0$. Applying (\ref{eq3.1}) to $\Omega_N$,
we obtain
$$
\frac{A|\Omega_N|}{M}\leq \sum_{j=1}^{q}D^{+}(\Lambda_j)
$$
and the result thus follows by taking $N\rightarrow \infty$.
\end{pf}

\medskip

Using the previous lemma, we now prove Theorem \ref{th0.1}.

\medskip

\noindent{\it Proof of Theorem \ref{th0.1}.} Suppose that there exists a frame of
 windowed exponentials $\bigcup_{j=1}^{q}{\mathcal E}(g_j,\Lambda_j)$ for $L^2(\Omega)$
with $|\Omega|=\infty$. By Lemma \ref{lem3.1}(ii), $D^{+}\Lambda_j=\infty$ for all $j$.
On the other hand, we consider $\Omega_N$ as in the proof of Lemma \ref{lem3.1}(ii)
with $N$ large enough so that $|\Omega_N|>0$. The windowed exponentials continue still
form a frame for $L^2(\Omega_N)$. By Lemma \ref{lem3.1}(i) applied to $\Omega_N$,
$D^{+}\Lambda_j<\infty$ for all $j$. This leads us to a contradiction and hence there
cannot be any frame of windowed exponentials for $L^2(\Omega)$ if $|\Omega|=\infty$.
\eproof

\medskip

From now on, we assume $|\Omega|<\infty$. The estimates obtained in following theorem
are the main technical tools used to
characterizating when  a system of windowed exponentials forms a Bessel sequence or a frame.
It further gives us explicit relationships between the frame bounds, the Beurling densities, and
the essential supremum and infimum of the moduli of the windows when dealing with a Bessel sequence
or frame of windowed exponentials.
Now, given a set  $\bigcup_{j=1}^{q}{\mathcal E}(g_j,\Lambda_j)$, we recall the
definition of $J$ and define a related index set $J'$:
\begin{equation}\label{eq3.5}
J = \{j: \|g_j\|_{\infty}<\infty\}, \ J' =J\cap \{j: D^{+}(\Lambda_j)>0\}
\end{equation}

\medskip

\begin{theorem}\label{th3.1}
Let $\Omega\subset {\Bbb R}^d$ such that $|\Omega|<\infty$ and
let $\bigcup_{j=1}^{q}{\mathcal E}(g_j,\Lambda_j)$ be a set of windowed exponentials in $L^2(\Omega)$.
\medskip

\noindent{\rm (i)} Let the collection $\bigcup_{i=1}^{q}{\mathcal E}(g_i,\Lambda_i)$
 form a Bessel sequence in $L^2(\Omega)$ with Bessel constant $B$ and suppose,
furthermore, that $D^{+}(\Lambda_j)>0$ for some $j\in\{1,\dots, q\}$.
Then, $|g_j|\leq \sqrt{B/D^{+}(\Lambda_j)}$ almost everywhere on $\Omega$.

\medskip

\noindent{\rm (ii)} If the collection  $\bigcup_{i=1}^{q}{\mathcal E}(g_i,\Lambda_i)$
 is a frame of
windowed exponentials in $L^2(\Omega)$ with frame bound $A,B (A<B)$, then we have
the inequalities
\begin{equation}\label{eq3.2}
\sqrt{A/D^{+}(\sum_{j\in J'}\delta_{\Lambda_j})}\leq \max_{j\in J'}|g_j|\leq
\max_{j\in J'}\sqrt{B/D^{+}(\Lambda_j)}
\end{equation}
almost everywhere on $\Omega$.
\end{theorem}

\medskip

\begin{pf}
(i) Suppose that  $\bigcup_{i=1}^{q}{\mathcal E}(g_i,\Lambda_i)$ is a Bessel
sequence in $L^2(\Omega)$, then ${\mathcal E}(g_i,\Lambda_i)$ are Bessel
sequences in $L^2(\Omega)$ for all $i$. Now if $D^{+}(\Lambda_j)>0$ for some $j$,
consider the measurable set
$$
E_{M} = \{x\in\Omega:M\leq |g_j(x)|\}.
$$
Note that  $|E_M|\leq \|g_j\|_{L^2(\Omega)}^2/M^2<\infty$, so $f := \chi_{E_{M}}\in
 L^2(\Omega)$ and $f$ satisfies the assumption in Proposition \ref{prop2.2}(i). Hence,
$$
 D^{+}\left((\int_{E_M}|g_j|^2)\cdot\delta_{\Lambda_j}\right)\leq B|E_M|.
$$
Since $\int_{E_M}|g_j|^2\geq M^2|E_M|$, we obtain that
$$
|E_M|\,(M^2 D^{+}(\Lambda_j)-B)\leq 0.
$$
If $M>\sqrt{B/D^{+}(\Lambda_j)}$, we would have $(M^2 D^{+}(\Lambda_j)-B)> 0$ which would
force $|E_M|$ to be zero. Hence,  $|g_j|\leq \sqrt{B/D^{+}(\Lambda_j)}$ a.e. on $\Omega$.
This establishes (i).

\medskip

\noindent{(ii)} From (i), we have shown for those $j\in J'$, $|g_j|$ is
essentially bounded above by $\sqrt{B/D^{+}(\Lambda_j)}$, from which the
second inequality in (\ref{eq3.2}) follows. It remains to establish the
first inequality in (\ref{eq3.2}). We consider, for $\epsilon>0$, the set
$$
F_{\epsilon} := \bigcap_{j\in J'}\{x\in\Omega: |g_j|< \epsilon\}.
$$
Define $f = \chi_{F_{\epsilon}}$. Then $\int_{\Omega}|fg_j|^2 \leq
\epsilon^2|F_{\epsilon}|\leq \epsilon^2|\Omega|<\infty$. By Proposition \ref{prop2.2},
\begin{equation}\label{eq3.3}
A|F_{\epsilon}|\leq D^{-}\left(\mu_f\right)\leq D^{+}\left(\mu_f\right)
\end{equation}
where $\mu_f=\sum_{j=1}^{q}\left(\int_{\Omega}|fg_j|^2\right)\,\delta_{\Lambda_j}$.
Note that if $g_j$ is not essentially bounded above, then $D^{+}(\Lambda_j)=0$ by (i) above.
If $j\in J\setminus J'$, then $D^{+}(\Lambda_j)=0$ also by the definition of $J'$.
We can now use  Proposition \ref{prop2.1} to conclude that
$$
D^{+}(\mu_f) = D^{+}\left(\sum_{j\in J'}(\int_{\Omega}|fg_j|^2)\, \delta_{\Lambda_j}\right).
$$
Note from the definition of $F_{\epsilon}$ that for those $j\in J'$,
$\int_{\Omega}|fg_j|^2\leq \epsilon^2 |F_{\epsilon}|$ and
$$
D^{+}(\mu_f)\leq \epsilon^2 |F_{\epsilon}|\, D^{+}(\sum_{j\in J'}\delta_{\Lambda_j}).
$$
Using (\ref{eq3.3}), we obtain that
$$
|F_{\epsilon}|\left(A- \epsilon^2D^{+}(\sum_{j\in J'}
\delta_{\Lambda_j})\right)\leq 0.
$$
This shows that $|F_{\epsilon}|=0$ and $\Omega\setminus F_{\epsilon}$ has
full measure in $\Omega$ if $\epsilon<\sqrt{A/D^{+}(\sum_{j\in J'}\delta_{\Lambda_j})}$.
As $\Omega\setminus F_{\epsilon_0} = \bigcap_{n=1}^{\infty}\Omega\setminus F_{\epsilon_0-1/n}$
and $|\Omega|<\infty$,  $\Omega\setminus F_{\epsilon_0}$ has full measure in $\Omega$
for $\epsilon_0 =\sqrt{A/D^{+}(\sum_{j\in J'}\delta_{\Lambda_j})}$. Note that
$$
\left\{x\in\Omega: \max_{j\in J'}|g_j(x)|\geq \epsilon_0\right\} = \bigcup_{j\in J'}
\left\{x\in\Omega:|g_j(x)|\geq \epsilon_0\right\} = \Omega\setminus F_{\epsilon_0}.
$$
This establishes the lower bound.
\end{pf}

\medskip

Now we prove Theorem \ref{th0.2}.

\medskip

\noindent{\it Proof of Theorem \ref{th0.2}.} Suppose that there exists a
frame of windowed exponentials $\bigcup_{j=1}^{q}{\mathcal E}(g_j,\Lambda_j)$
for $L^2(\Omega)$. By Lemma \ref{lem3.1}(i), $D^{+}(\Lambda_j)<\infty$ for all $j$
and $D^{+}(\Lambda_j)>0$ for at least one such $j$. Using Theorem \ref{th3.1}(i)
the corresponding $g_j$ is essentially bounded above on $\Omega$. Hence, $J'$ and
therefore $J$ in (\ref{eq3.5}) are non-empty. By Theorem \ref{th3.1} (ii),
$\max_{j\in J'} |g_j|$ is essentially  bounded away from 0 on $\Omega$ and,
since $J\supset J'$, so is $\max_{j\in J} |g_j|$.  This shows
the necessity of that condition.

\medskip

Suppose now  $m\leq \max_{j\in J}|g_j|$. Since $J$ is a finite set, the
definition of $J$ shows that  $\max_{j\in J}|g_j|\leq M$ for some $M<\infty$.
As $\Omega$ is bounded, we can cover $\Omega$ by a cube $Q_R$.
We know that $L^2(Q_R)$ has a Fourier frame (in fact an orthonormal basis),
which we denote by $\{e^{2\pi i \langle\lambda,\cdot\rangle}\}_{\lambda\in\Lambda}$.
Define
$$
\Lambda_{j} =\left\{
               \begin{array}{ll}
                 \Lambda, & \hbox{$j\in J$;} \\
                 \{0\}, & \hbox{$j\not\in J$.}
               \end{array}
             \right.
$$
To prove the upper bound in the frame inequality, we note that, for $j\in J$,
$$
\int_{\Omega}|fg_j|^2\leq M^2\int_{\Omega}|f|^2<\infty,\quad f\in L^2(\Omega).
$$
Hence, denoting by $B$ the Bessel constant of the sequence of the exponentials associated
with $\Lambda$ for $L^2(Q_R)$, we have
$$
\sum_{j\in J}\sum_{\lambda\in\Lambda}\left|\int_{\Omega}f(x)\overline{g_j(x)}
e^{-2\pi i \langle\lambda,x\rangle}dx\right|^2\leq \sum_{j\in J} B\|fg_j\|_{L^2(Q_R)}^2
\leq B\,q\,M^2\,\|f\|_{L^2(\Omega)}^2
$$
for all $f\in L^2(\Omega)$ (we take $f =0$ on $Q_R\setminus \Omega$).
While for $j\not\in J$, we simply use Cauchy-Schwarz inequality to obtain
$$
\sum_{j\not\in J}|\int_{\Omega}f(x)\overline{g_j(x)}dx|^2\leq q|\Omega|
\max_j\{\|g_j\|_{L^2(\Omega)}^2\}\,\|f\|_{L^2(\Omega)}^2.
$$
Hence, combining these  last
two inequalities, we obtain
$$
\sum_{j=1}^{q}\sum_{\lambda\in\Lambda}\left|\int_{\Omega}f(x)
\overline{g_j(x)}e^{-2\pi i \langle\lambda,x\rangle}dx\right|^2\leq q(BM^2+
|\Omega|\max_j\{\|g_j\|_{L^2(\Omega)}^2\})\,\|f\|_{L^2(\Omega)}^2.
$$
which yields the upper bound in the frame inequality

\medskip

To establish the lower bound, we note that we can
remove those $\Lambda_j$ with $j\not\in J$ in the sum appearing in the middle of the
ineqaulities in (\ref{eq0.1}).
Now, using the fact that $m\leq \max_{j\in J}|g_j|$ a.e., we deduce that
$$
\left|\Omega\setminus\bigcup_{j\in J}\{|g_j|\geq m\} \right| =0.
$$
We can replace, if necessary, the sets $\{|g_j|\geq m \}$ with $j\in J$
by subsets $T_j$, $j\in J$, which still cover $\Omega$ and are pairwise disjoint.
Denoting by $A$ lower frame bound for the set of exponentials on $L^2(Q_R)$
with associated frequencies in $\Lambda$, we have thus
$$
\begin{aligned}
\sum_{j\in J}\sum_{\lambda\in\Lambda}\left|\int_{\Omega}f(x)\overline{g_j(x)}
e^{-2\pi i \langle\lambda,x\rangle}dx\right|^2\geq& A\sum_{j\in J}\|fg_j\|_{L^2(\Omega)}^2\\
\geq& A\sum_{j\in J}\int_{T_j}|fg_j|^2\\
\geq& A m^2 \sum_{j\in J}\int_{T_j}|f|^2= A m^2 \int_{\Omega}|f|^2,
\end{aligned}
$$
where the  pairwise disjointness of the sets $T_j$, $j\in J$, and their covering property
is used to obtain the last equality. This yields the lower bound and proves our claim.

\eproof

\medskip

\begin{Example}
  {\rm On the interval $[0,1]$, let $g_1(x) = x^{\alpha}$ and $g_2(x) =
 (1-x)^{\alpha}$ with $\alpha\geq 0$. Then }
$$
\max\{g_1(x), g_2(x)\} =
\left\{
  \begin{array}{ll}
                             (1-x)^{\alpha}, & \hbox{$0\leq x\leq1/2$;} \\
                             x^{\alpha}, & \hbox{$1/2<x\leq 1$.}
  \end{array}
\right.
 \   \mbox{if} \  \alpha\geq 0,
$$
{\rm Hence, $1/2\leq\max\{g_1,g_2\}\leq 1$ on $[0,1]$. We can produce a frame of exponentials
on $[0,1]$ by taking $\Lambda_1 = \Lambda_2 = {\Bbb Z}$ for instance.}

\medskip

{\rm  On the other hand, the functions $g_3(x) = x^{\beta}$ and $g_4(x) = (1-x)^{\beta}$  are in $L^2([0,1])$
if $-1/2<\beta<0$. Since they are both unbounded,
they cannot be used to produce windowed exponentials for $L^2([0,1])$. }

\medskip

{\rm If we now consider the collection $\{g_1,g_2,g_3,g_4\}$. We note that the
sub-collection $\{g_1,g_2\}$  satisfies  Theorem \ref{th0.2}, so we can use this collection to
form a frame of windowed exponentials and the unbounded functions $g_3$ and $g_4$ are redundant.}
\end{Example}

We end this section with a remark about generalized frames of windowed exponentials.

\begin{Rem}\label{rem4.1}
{\rm Let $g_1,\cdots,g_q\in L^2(\Omega)\setminus\{0\}$, and $\mu_1,\cdots,\mu_q$ be
locally finite Borel measures on ${\Bbb R}^d$, we say that the collections $\bigcup_{i=1}^{q}{\mathcal E}(g_j,\mu_j)$
form a {\it generalized frame of windowed exponentials} for $L^2(\Omega)$ if we can
find $0<A,B<\infty$ such that}
$$
A\|f\|^2_{L^2(\Omega)}\leq\sum_{j=1}^{q}\int\left|\int_{\Omega} f(x)
\overline{g_j(x)}e^{-2\pi i \langle\lambda,x\rangle}dx\right|^2d\mu_{j}(\lambda)\leq B\|f\|^2_{L^2(\Omega)}
$$
{\rm  for all $f\in L^2(\Omega)$. Similar generalized frames were also studied in \cite{[DHW]}.
In the case where $\mu_j = \delta_{\Lambda_j}$, $j=1,\cdots, q$, we recover the system of the
windowed exponentials defined in the introduction. As Theorem \ref{th2.3} and
Proposition \ref{Prop2.2-} are true for general measures $\mu$, all the arguments in
 Proposition \ref{prop2.2} and this section  holds by directly replacing $\delta_{\Lambda_j}$
with $\mu_j$. Therefore,    Theorem \ref{th0.1} and \ref{th0.2} actually holds for generalized  frames of
windowed exponentials.}
\end{Rem}

\bigskip

\section{unbounded sets of finite measures}
 We know that a bounded set in ${\Bbb R}^d$ can be covered by a hypercube and, in particular,
the  orthonormal bases of exponentials defined on the cube that we mentioned earlier
 will generate a tight frame for that set when restricted to it.
When the set is unbounded but is of finite Lebesgue measure, such argument generally
fails unless the set considered has some special properties such as in the following proposition.
The result is known (see e.g. \cite{[GaL]} for the one-dimensional case),
but we provide here a simple proof for the reader's convenience.

 \medskip

\begin{Prop}\label{prop5.1}
Let  $\Omega$  be a set of finite Lebesgue measure (bounded or unbounded).
Let $\Gamma$ be a lattice in ${\Bbb R}^d$ with
$\Gamma^{\ast} = \{\lambda: \langle\lambda,\gamma\rangle\in{\Bbb Z}\}$ being
its dual lattice. Then, the following are equivalent.

\medskip

(i) $\sum_{\gamma\in\Gamma}\chi_{\Omega}(x+\gamma)\leq 1$ almost everywhere on ${\Bbb R}^d$.

\medskip

(ii) The collection $\{e^{2\pi i \langle\lambda,\cdot\rangle}\}_{\lambda\in\Gamma^{\ast}}$
is a (tight) Fourier frame for
$L^2(\Omega)$.
\end{Prop}

\begin{pf}
(i) $\Longrightarrow$ (ii). Let $Q$ be a fundamental domain of $\Gamma$ and let $f\in L^2(\Omega)$. We have
 \begin{equation}\label{eq5.1}
\int_{\Omega}\,|f(x)|^2\,dx =
\int_{Q}\,\sum_{\gamma\in\Gamma}\,\chi_{\Omega}(x+\gamma)\,|f(x+\gamma)|^2\,dx.
\end{equation}
From the assumption $\sum_{\gamma\in\Gamma}\chi_{\Omega}(x+\gamma)\leq1$
almost everywhere, we have that almost every $x\in Q$, there exists at
most one $\gamma_x\in\Gamma$ such that $\chi_{\Omega}(x+\gamma_x)=1$ which implies that
 $$
\sum_{\gamma\in\Gamma}\chi_{\Omega}(x+\gamma)|f(x+\gamma)|^2=
\big|\sum_{\gamma\in\Gamma}\chi_{\Omega}(x+\gamma)f(x+\gamma)\big|^2 \ \mbox{a.e. on $Q$}.
$$
It is well known that the system
$\{e^{2\pi i \langle\lambda,\cdot\rangle}\}_{\lambda\in\Gamma^{\ast}}$
 forms an orthogonal basis for $L^2(Q)$ (\cite{[Fu]}). Combining this fact with (\ref{eq5.1}),
we obtain
\begin{align*}
\int_{\Omega}|f(x)|^2\,dx &= \int_{Q}\,\big|\sum_{\gamma\in\Gamma}\,\chi_{\Omega}(x+\gamma)\,
f(x+\gamma)\big|^2\,dx \\
&= |Q| \sum_{\lambda\in\Gamma^{\ast}}\left|\int_{Q}
\sum_{\gamma\in\Gamma}\chi_{\Omega}(x+\gamma)f(x+\gamma)e^{-2\pi i \langle\lambda,x\rangle}dx\right|^2.
\end{align*}
Since $\langle\lambda,\gamma\rangle\in{\Bbb Z}$ for $\lambda\in{\Gamma}^{\ast}$,
it follows directly that
\begin{align*}
&\int_{Q}\,\sum_{\gamma\in\Gamma}\,
\chi_{\Omega}(x+\gamma)f(x+\gamma)\,e^{-2\pi i \langle\lambda,x\rangle}\,dx\\
&=e^{2\pi i \langle\lambda,\gamma\rangle}\int_{\Omega}f(x)e^{-2\pi i \langle\lambda,x\rangle}\,dx
=\int_{\Omega}\,f(x)\,e^{-2\pi i \langle\lambda,x\rangle}\,dx
\end{align*}
 which implies
$$
\int_{\Omega}|f(x)|^2\,dx =  |Q| \sum_{\lambda\in\Gamma^{\ast}}
\left|\int_{\Omega}\,f(x)\,e^{-2\pi i \langle\lambda,x\rangle}\,dx\right|^2.
$$

\medskip

(ii) $\Longrightarrow$ (i). Let $Q$ be a fundamental domain of $\Gamma$. Proving the statement in (i)
is equivalent to showing that
$$
\sum_{\gamma\in\Gamma}\chi_{\Omega}(x+\gamma)\leq 1 \quad \text{a.e.~on}\,\,Q
$$
since the term on the left-hand side of the inequality above is $\Gamma$-periodic.
This is in turn
equivalent to showing
that
$$
|(\Omega-\gamma) \cap (\Omega-\gamma')\cap Q|=0\quad \text{for
all}\,\, \gamma,\gamma'\in\Gamma\,\,\text{with}\,\, \gamma\neq\gamma'.
$$
We argue by contradiction and  suppose that   there exists $\gamma,\gamma'\in\Gamma$
such that $\gamma\neq\gamma'$ and
$|(\Omega-\gamma) \cap (\Omega-\gamma')\cap Q|>0$.
Let $E =(\Omega-\gamma) \cap (\Omega-\gamma')\cap Q$ and consider
$f =  \chi_{E+\gamma}-\chi_{E+\gamma'}\in L^2(\Omega)$. Note that $f$ is a
 non-zero function in $L^2(\Omega)$ since $|E|>0$. On the other hand, for
all $\lambda\in\Gamma^{\ast}$,
$$
\int_{\Omega}f(x)e^{-2\pi i \langle\lambda,x\rangle}dx = \int_{E+\gamma}
e^{-2\pi i \langle\lambda,x\rangle}dx-e^{2\pi i \langle\lambda,\gamma-\gamma'\rangle}
\int_{E+\gamma} e^{-2\pi i \langle\lambda,x\rangle}dx=0.
$$
This shows the system
$\{e^{2\pi i \langle\lambda,\cdot\rangle}\}_{\lambda\in\Gamma^{\ast}}$ is
incomplete in $L^2(\Omega)$ and hence cannot be a Fourier frame for $L^2(\Omega)$.
This contradicts
the assumption in (ii).
\end{pf}

\medskip

We will now proceed to prove Theorem \ref{th0.3}. To this end, we need another approach using certain
convolution identities due to Kolountzakis \cite{[K1],[K2]}. If a function
$f\in L^1({\Bbb R}^d)$ and a countable set $\Lambda\subset{\Bbb R}^d$ are such that
$$
\delta_{\Lambda}\ast f(x) = \sum_{\lambda\in\Lambda}f(x-\lambda)
=w \ \mbox{ a.e. $x\in{\Bbb R}^d$},
$$
$f$ is said to {\it tile by $\Lambda$ at level $w$} for some constant $w$
(see \cite{[K1],[K2]}).
Kolountzakis \cite[Theorem 2]{[K1]} proved the following result in which
 $\widehat{\delta_{\Lambda}}$  denotes the distributional Fourier transform
of $\delta_{\Lambda}$ (as a tempered distribution).

 \begin{Prop}\label{prop5.2}
 Suppose that $f\in L^1({\Bbb R}^d)$ satisfies $f\geq 0$, $\widehat{f}\geq 0$ and assume that
the support of $\widehat{f}$ is compact. Then, if $f$ tiles by $\Lambda$ at some level $w>0$, we have
\begin{equation}\label{eq5.2}
\mbox{supp} \ \widehat{\delta_{\Lambda}}\subset \{x:\widehat{f}(x)=0\}\cup\{0\}.
\end{equation}
 \end{Prop}

%By [K2, Theorem 1.11], it is in general not defined by a measure, but it is locally a measure in the following sense: If $\Lambda$ is a discrete set with a uniform density (i.e. $D^{\pm}\Lambda = \rho>0$), then for $a\in {\Bbb R}^d$ and $R>0$ such that
%$$
%\mbox{supp}\widehat{\delta_{\Lambda}}\cap B_{R}(a) = \{a\},
%$$
%($B_R(a)$ is the ball of radius $R$ centered at $a$), $\widehat{\delta_{\Lambda}} = c\delta_a$, for some constant $c$ on $B_R(a)$.

\medskip

\noindent{\it Proof of Theorem \ref{th0.3}.} We argue by contradiction.
Suppose that we are given a set $\Omega$ of finite measure and $R>0$ such that
$|\Omega\cap\Omega+x|>0$ for all $|x|>R$ and that $L^2(\Omega)$ admits a tight
frame $\{e^{2\pi i \langle\lambda,\cdot\rangle}\}_{\lambda\in\Lambda}$ with
frame constant $A$. Consider $\Omega_N = \Omega\cap [-N,N)^d$. Clearly, when restricted to $\Omega_N$,
this tight frame produces a tight frame  for $L^2(\Omega_N)$. Hence,
applying the definition of tight frame to the function
$\chi_{\Omega_N}e^{2\pi i \langle\xi,\cdot\rangle}$, $\xi \in {\Bbb R}^d$, we obtain
$$
\left(\delta_{\Lambda}\ast |\widehat{\chi_{\Omega_N}}|^2\right)(\xi) =
\sum_{\lambda\in\Lambda}\,|\widehat{\chi_{\Omega_N}}(\xi-\lambda)|^2= A\,|\Omega_N|,\quad \xi \in  {\Bbb R}^d.
$$
Letting $f_N = |\widehat{\chi_{\Omega_N}}|^2$,  we have
$\widehat{f_N} = \chi_{\Omega_N}\ast\widetilde{\chi_{\Omega_N}}$
where $\widetilde{\chi_{\Omega_N}}(x) = \chi_{\Omega_N}(-x)$.
A simple calculation shows that
 $$
\widehat{f_N}(x) = |\Omega_N\cap \Omega_N+x|,\quad x\in {\Bbb R}^d.
$$
Since $\Omega_N$ is bounded,  $\widehat{f_N}$ has compact support
also and $\widehat{f_N}\ge 0$. By Proposition \ref{prop5.2},
$$
\mbox{supp} \ \widehat{\delta_{\Lambda}}\subset \{x:|\Omega_N\cap \Omega_N+x|=0\}
\cup\{0\} \ \mbox{for all $N>0$}.
$$
Note that, since $\Omega_N$ is an increasing sequence of sets
whose union is $\Omega$, $|\Omega_N\cap\Omega_N+x|$ converges pointwise to
$|\Omega\cap\Omega+x|$ as $N\rightarrow\infty$. Hence, using our assumption on $\Omega$,
for all $x$ such that $|x|>R$, there exists $N$ such that
$|\Omega_N\cap\Omega_N+x|>0$. This means that
$\bigcap_N\{x:|\Omega_N\cap \Omega_N+x|=0\}\subset \{x\in{\Bbb R}^d: |x|<R\}$.
Therefore,
$$
\mbox{supp}\ \widehat{\delta_{\Lambda}}\subset
\bigcap_{N=1}^{\infty}\{x:|\Omega_N\cap \Omega_N+x|=0\}\cup\{0\}
\subset \{x\in{\Bbb R}^d:|x|\leq R\}
$$
showing that the support of $\widehat{\delta_{\Lambda}}$ is compact.
This leads to a contradiction since, by the
Paley-Wiener-Schwartz theorem (\cite[p.199]{[R]}), any tempered distribution
whose Fourier transform is compactly supported must be the restriction to
${\Bbb R}^d$ of an entire analytic function, but here $\delta_{\Lambda}$
is a purely discrete measure on ${\Bbb R}^d$.
\eproof

\medskip

 We can strengthen Theorem \ref{th0.3} to a more general setting.
We say that $L^2(\Omega)$ admits a {\it generalized tight frame of exponentials} if
there exists a locally finite Borel measure $\mu$ on ${\Bbb R}^d$ and
a constant $A>0$ such that
\begin{equation}\label{eq4.2}
\int_{{\Bbb R}^d}|\widehat{f}(\lambda)|^2d\mu(\lambda) =
A\,\int_{\Omega}|f(x)|^2dx,  \ \ f\in L^2(\Omega).
\end{equation}
In this situation, $\mu$ is called a {\it tight frame measure}
for $L^2(\Omega)$ (see \cite{[DHW],[DL]}). It is clear that if
$\mu$ is the Lebesgue measure on ${\Bbb R}^d$, then, by the
Plancherel theorem, $\mu$ is a tight frame measure for $L^2(\Omega)$
for any $\Omega$. In Theorem \ref{th5.1}, we give a necessary and
sufficient condition for the Lebesgue measure on ${\Bbb R}^d$ to
be the only tight frame measure for $L^2(\Omega)$. An analogous problem
was also considered by the first named author in the setting of
Gabor analysis (\cite{[Ga0]}). In particular, we need a lemma in (\cite[Lemma 4.5]{[Ga0]}).

\begin{Lem}\label{lem5.1}
Let $\mu$ be a positive translation-bounded  measure on ${\Bbb R}^d$.
Suppose that for some $r>0$ and some $\tau\in{\Bbb R}^d$,
$$
\mbox{supp} \ \widehat{\mu}\cap B_{r}(\tau) = \{\tau\},
$$
where $B_{r}(\tau)$ is the ball of radius $r$ centered at $\tau$.
Then, there exists $a\in {\Bbb C}$ such that
$$
\widehat{\mu} = a\delta_{\tau} \ \mbox{on} \ B_{r}(\tau).
$$
\end{Lem}

\medskip

\begin{Rem}\label{rem5.1+}
{\rm If $\tau=0$ in the previous lemma, then $a>0$. Indeed,
if $\varphi\in C^{\infty}_{0}({\Bbb R}^d)$ is supported on
$B_{r/2}(0)$, then  the support of $\varphi\ast\widetilde{\varphi}$
(recall that $\widetilde{\varphi}(x) = \varphi(-x)$) is contained
in $B_r(0)$ and $a\|\varphi\|_2^2=\langle\widehat{\mu},
\varphi\ast\widetilde{\varphi}\rangle = \int|\widehat{\varphi}(\xi)|^2d\mu(\xi)>0$. Hence, $a>0$.}
\end{Rem}

\medskip

We also need to use Proposition \ref{prop5.2} with $\delta_{\Lambda}$
 replaced by $\mu$. This is possible by a simple modification of the
argument in the proof given in \cite[Theorem 2]{[K1]}. We leave the details
to the interested reader.

%since we know from the proof that the assumptions on $f$ in Proposition \ref{prop5.2} are to guarantee $\varphi/\widehat{f}\in L^1({\Bbb R}^d)$ for any $\varphi\in C^{\infty}_c(K^c)$, where $K$ denotes the set in the right hand side of (\ref{eq5.2}).  $\widehat{\mu}(\varphi) =0$ can then be justified using Fubini theorem by an obvious modification.

\medskip

\begin{theorem}\label{th5.1}
Let $\Omega$ be a measurable subset of ${\Bbb R}^d$ with $|\Omega|<\infty$.
Then the following are equivalent.

(i) The only tight frame measure for $L^2(\Omega)$ is the
Lebesgue measure on ${\Bbb R}^d$, up to a positive constant multiple.

(ii) $|\Omega\cap(\Omega+x)|>0$ for all $x\in{\Bbb R}^d\setminus\{0\}$.
\end{theorem}

\begin{pf}
(ii) $\Longrightarrow$ (i).  Suppose that $\Omega$ is a set of finite
measure satisfying (ii) and that $L^2(\Omega)$ admits a tight
frame measure. Consider $\Omega_N = \Omega\cap [-N,N)^d$. By restriction,
this tight frame measure continues to be a tight frame measure for $L^2(\Omega_N)$.
Replacing $f$ by $\chi_{\Omega_N}\,e^{2\pi i \langle\xi,\cdot\rangle}$, $\xi \in {\Bbb R}^d$,  in (\ref{eq4.2}), we obtain
$$
\left(\mu\ast|\widehat{\chi_{\Omega_N}}|^2\right)(\xi) = A|\Omega_N|,\quad \xi \in {\Bbb R}^d.
$$
 Using a similar argument as in the proof of Theorem \ref{th0.3} with
$\delta_{\Lambda}$ replaced by $\mu$, we  deduce that
 $\mbox{supp}\ \widehat{\mu} = \{0\}$ since $\bigcap_{N=1}^{\infty}
\{x:|\Omega_N\cap \Omega_N+x|=0\}=\{0\}$. By Lemma \ref{lem5.1} and
the remark following it, it follows that $\widehat{\mu} = a\,\delta_0$.
This shows that the only tight frame measures for $L^2(\Omega)$
are the positive multiples of the Lebesgue measure on ${\Bbb R}^d$.

\medskip

(i) $\Longrightarrow$ (ii). Suppose that there exists $x_0\neq 0$
such that $|\Omega\cap(\Omega+x_0)| =0$. Then it is easy to see
that $|\Omega\cap(\Omega-x_0)| =0$ also. Define the following locally finite measure
$$
d\mu(\xi) = \left(1+\frac{e^{2\pi i \langle x_0,\xi\rangle}+
e^{-2\pi i \langle x_0,\xi\rangle}}{2}\right)\,d\xi.
$$
It is a positive measure since the density is equal to $1+
\cos(2\pi  \langle x_0,\xi\rangle)\geq 0$. We now claim that
it is another tight frame measure for $L^2(\Omega)$ and this will prove our claim.
Indeed, for any $f\in L^2(\Omega)$,
$$
\int|\widehat{f}(\xi)|^2e^{2\pi i \langle x_0,x\rangle}d\xi =
\int\widehat{f}(\xi)e^{2\pi i \langle x_0,x\rangle}\overline{\widehat{f}(\xi)}d\xi
=\int f(x+x_0)f(x)dx.
$$
As $f$ is supported on $\Omega$, the integrand $f(x+x_0)f(x)$
is supported on the intersection $\Omega\cap (\Omega-x_0)$ which has zero Lebesgue measure.
This shows the integral above is $0$. The same also applies to  $\int|\widehat{f}(\xi)|^2e^{-2\pi i \langle x_0,x\rangle}d\xi$. This shows that
$$
\int|\widehat{f}(\xi)|^2d\mu(\xi) = \int|\widehat{f}(\xi)|^2d\xi = \int_{\Omega}|f(x)|^2dx,
$$
which completes the proof.
 \end{pf}
 %To conclude the proof, we know from the paragraph before this proof that $\widehat{\delta_{\Lambda}} = c\delta_0$ is a measure on all $B_R(0)$. Clearly, $c$ is independent of $R$. Thus, $\widehat{\delta_{\Lambda}}$ is the Dirac measure at $0$ of some constant $c$. Hence, $\delta_{\Lambda}$ can only be the Lebesgue measure on ${\Bbb R}^d$. This is a contradiction since $\delta_{\Lambda}$ is a purely discrete measure.

\medskip

We now give an example of a set $\Omega$ of finite measure such
that any two translates of $\Omega$ always intersect on a set of
 positive measure. We only present an example on ${\Bbb R}$ as higher dimensional example can easily
be constructed from it.

\begin{Example}
The set $\Omega = \bigcup_{n\in{\Bbb Z}}\left([-\frac{1}{2^{|n|}},\frac{1}{2^{|n|}}]+n\right)$
 is a set of finite Lebesgue measure satisfying $|\Omega\cap(\Omega+x)|>0$ for all $x\in{\Bbb R}$.

If we let $\Omega_k = [-1,1]\cup\bigcup_{|n|>k}\left([-\frac{1}{2^{|n|}},
\frac{1}{2^{|n|}}]+n\right)$ for $k\geq 4$, then $\Omega_k$ has
finite measure and $|\Omega_k\cap(\Omega_k+x)|>0$ for all $|x|\geq k$,
but the set of $x$ such that $\Omega_k\cap(\Omega_k+x)=\emptyset$
has positive measure.
\end{Example}
\begin{pf}
The finiteness of the Lebesgue measure of $\Omega$ and $\Omega_k$
are clear. Let $x\in{\Bbb R}$ and let $n$ be the unique
integer such that $n\leq x<n+1$, then $\Omega+x\supset[x,x+1]$ and
the interval $[n+1-\frac{1}{2^{|n+1|}},n+1+\frac{1}{2^{|n+1|}}]$
intersects the interval $[x,x+1]$ on a set of positive measure.
This shows that $\Omega\cap\Omega+x$ has positive Lebesgue measure.

\medskip

Using the same method above for $\Omega_k$, we can show
$|\Omega_k\cap(\Omega_k+x)|>0$ for all $|x|\geq k$.
Now consider $x=5/2$. Then $(\Omega_k+x) =
[3/2,7/2]\cup\bigcup_{|n|>k}\left([-\frac{1}{2^{|n|}},\frac{1}{2^{|n|}}]+n+5/2\right)$.
As $[3/2,7/2]$ does not intersect $\Omega_k$ if $k>4$ and the
lengths of the remaining intervals centered at $n+1/2$ are all less
than $1/2^4$, $\Omega_k+5/2$ is  disjoint from all the intervals in $\Omega_k$.
Moreover, $\Omega_k$ and $\Omega_k+5/2$ are at a positive distance from each other.
Therefore, for all $x$ close to $5/2$, $\Omega_k$ and $\Omega_k+x$ also has this property.
Hence, this shows the set of $x$ such that $\Omega_k\cap(\Omega_k+x)=\emptyset$ has positive measure.
\end{pf}

\medskip

We conclude this section with some remarks.

\begin{Rem}
{\rm  (1). It is unknown whether (non-tight) Fourier frames always exist
for unbounded sets of finite measure.  This problem was addressed
earlier in \cite{[OU]}. From all the approaches we tried in which
we assume that any two translates intersect on a set of positive
measure, we cannot formulate a definite conjecture to this problem. For instance,}

{\rm (i) It is even possible construct sets of finite measure such
that any finite number of translates intersect with positive measure.
For these sets, we can show that if a Fourier frame exists with some
frequency set $\Lambda$, then the set $\Lambda \ (\mbox{mod} \  \Gamma)$
has to be dense in the fundamental domain of any lattice $\Gamma$.}

{\rm (ii) On the other hand, Matei and Meyer \cite{[MM]} recently
constructed from simple quasicrystals a universal Fourier frame. This
means that the frequency set $\Lambda$ will form a Fourier frame
on any $L^2(K)$ such that  $D^{-}\Lambda>|K|$ and $K$ is compact
with boundary measure 0. Their method may be extendable to cover
our sets.}

\medskip

\noindent{\rm (2). Another problem of a similar nature asks whether or not a
Fourier frame exists on the singular one-third Cantor measure.
In the existing methods, the construction of a Fourier frame is based
on the existence of a singular measure for which there exists an
orthonormal basis of exponentials  (\cite{[HLL],[DL]}). While it
is known that the one-third Cantor measure cannot admit any
exponential orthogonal basis (\cite{[JP]}), we are interested in
the existence of Fourier frames for a measure which genuinely
cannot be derived from some already existing tight frames.}
\end{Rem}

\bigskip

\section{Some Applications}
We now give some application of our result to other well-known types of
frames.

\medskip

{\it (I) Frame of translates}

\vspace{0.2cm}

Given $g_1,\cdots, g_m\in L^2({\Bbb R}^d)$ and associated  countable sets ${\mathcal J}_i$ in ${\Bbb R}^d$, $i=1,\cdots,m$,
it was shown in  \cite{[CDH]}
that $\bigcup_{j=1}^{m}\{g_j(x-t): t\in{\mathcal J}_j\}$ cannot be a
frame for $L^2({\Bbb R}^d)$. We give a simple proof of this fact, based on our previous results
concerning windowed exponentials.

\begin{theorem}
There is no frame  of the form $\bigcup_{j=1}^{m}\{g_j(x-t): t\in{\mathcal J}_j\}$ on $L^2({\Bbb R}^d)$.
\end{theorem}

\begin{pf}
Suppose there is a frame of the given form. Then the Fourier
transforms of the functions in the system will also form a frame for $L^2({\Bbb R}^d)$.
Since $\widehat{g_j(\cdot-t)}(\xi) =
\widehat{g_j}(\xi)\,e^{2\pi i \langle t,\xi\rangle}$,
this new system will be  in the form of windowed exponentials on $L^2({\Bbb R}^d)$, generated by a finite number
of windows, contradicting  Theorem \ref{th0.1}.
\end{pf}

\medskip

{\it (II) Frames of absolutely continuous measures}

\vspace{0.2cm}

Let $\mu$ be an absolutely continuous measures with compact support.
We write $d\mu(x) = \varphi(x)dx$, where $\varphi$ is its
Radon-Nikodym derivative. In \cite{[Lai]},
the second named author completely characterized the kind of
density such that the measure admits a Fourier frame. We say
that a measure $\mu$ has an associated frame of windowed exponentials
if we can find $\bigcup_{j=1}^{q}{\mathcal E}(g_j,\Lambda_j)$ with
 $g_j\in L^2(\mu)$ which forms a frame for $L^2(\mu)$. i.e.
\begin{equation}\label{eq4.1}
A\|f\|_{L^2(\mu)}^2\leq \sum_{j=1}^{q}\sum_{\lambda\in\Lambda_j}
\left|\int f(x)g_j(x)e^{-2\pi i \langle\lambda,x\rangle}d\mu(x)\right|^2
\leq B\|f\|_{L^2(\mu)}^2, f\in L^2(\mu).
\end{equation}
The following proposition shows that for an absolutely continuous
measure, the notion of frame of windowed exponentials associated
with the measure and the  frame of exponentials on the support of
the measure are equivalent.

\medskip

\begin{Prop}\label{prop4.1}
Let $\mu = \varphi(x)dx$ be an absolutely continuous measures
and let $\Omega = \{\varphi\neq 0\}$. Then $\bigcup_{j=1}^{q}{\mathcal E}(g_j,\Lambda_j)$
is a frame of exponentials for $L^2(\varphi dx)$ if and only if
$\bigcup_{j=1}^{q}{\mathcal E}(g_j\sqrt{\varphi},\Lambda_j)$ is a
frame of exponentials for $L^2(\Omega)$.
\end{Prop}

\begin{pf}
Suppose $\bigcup_{j=1}^{q}{\mathcal E}(g_j,\Lambda_j)$ is a
frame of exponentials of $L^2(\varphi dx)$, then for any $f\in L^2(\Omega)$,
we have  $\int_{\Omega}|\frac{f(x)}{\sqrt{\varphi(x)}}|^2\varphi(x) dx =
\int_{\Omega}|f(x)|^2dx<\infty$. Hence, we can replace f by $f/\sqrt{\varphi}$
in (\ref{eq4.1}), we obtain (\ref{eq0.1}).

\medskip

Conversely, if $\bigcup_{j=1}^{q}{\mathcal E}(g_j\sqrt{\varphi},\Lambda_j)$
is a frame of exponentials of $L^2(\Omega)$. Then for any $f\in L^2(\varphi dx)$,
we have $\int_{\Omega}|f\sqrt{\varphi}|^2 = \int_{\Omega}|f|^2\varphi dx<\infty$.
Therefore, replacing $f$ by $f\sqrt{\varphi}$ in (\ref{eq0.1}) and the
windows $g_j$ by $g_j\sqrt{\phi}$, we obtain (\ref{eq4.1}), which proves our claim.
\end{pf}

This leads to the following characterization for the frame of
windowed exponentials in $L^2(\varphi dx)$. The proof follows
easily from Theorem \ref{th0.2} and Proposition \ref{prop4.1}.

 \begin{theorem}\label{th4.2}
Let $\mu = \varphi(x)dx$ be an absolutely continuous measures
with $\Omega = \{\varphi\neq 0\}$ and let $g_j$, $j=1,2\cdots, q$ be a
finite set of functions in $L^2(\varphi dx)$. Then there exists
$\Lambda_j$ such that  $\bigcup_{j=1}^{q}{\mathcal E}(g_j,\Lambda_j)$
form a frame in $L^2(\varphi dx)$ if and only if there is a
sub-collection of functions $\{g_j\}_{j\in J}$, $J\subset \{1,\cdots q\}$
and constants $m,M$ with $0<m\leq M<\infty$ such that
$$
\frac{m}{\sqrt{\varphi}}\leq \max_{j\in J}|g_j|\leq \frac{M}{\sqrt{\varphi}}
$$
almost everywhere on $\Omega$.
\end{theorem}

If there is only one window $g = \chi_{\Omega}$ on $L^2(\varphi(x)dx)$,
Theorem \ref{th4.2} states that $\varphi$ must be bounded above and
bounded away from 0 on $\Omega$, which recovers the result in \cite{[Lai]}.

\medskip

{\it (III) Gabor frames}

\vspace{0.2cm}

Let $g\in L^2({\Bbb R}^d)$ and consider the Gabor system with lattice time-frequency shifts defined as follows.
$$
{\mathcal G}(g,a,b) = \{e^{2\pi i mb x} g(x-na):m,n\in{\Bbb Z}^d\}.
$$
It is well-known that if ${\mathcal G}(g,a,b)$ forms a frame for $L^2({\Bbb R}^d)$,
then $ab\leq 1$. The converse is in general false and characterizing the kind of
functions which form a Gabor frame is an important question. Since rescaling the
function $g$ does not affect the frame property, one can assume $b=1$ and $a\leq 1$.
One of the major tools in the theory of Gabor frames is the Zak transform. It is a unitary
mapping from $L^2({\Bbb R}^d)$ to $L^2([0,1]^{2d})$ defined by
$$
Zf(x,t) = \sum_{k\in{\Bbb Z}^d}f(x-k)e^{2\pi i \langle k,t\rangle}.
$$
If the previous definition of $Zf$ is extended to all of ${\Bbb R}^{2d}$, $Zf$ is
{\it quasiperiodic}  in the following sense:
$$
Zf(x,t+n) = Zf(x,t), \ Zf(x+n,t) = e^{2\pi i \langle n,t\rangle}Zf(x,t), \ \mbox{ $\forall$ $n\in{\Bbb Z}^d$.}
$$
It is also well known that if $a=1$, then ${\mathcal G}(g,1,1)$ is a Gabor
frame if and only if $0<A\leq |Zg|\leq B<\infty$ almost everywhere on $[0,1]^{2d}$
(see [G, p.157]). The following theorem  is a particular case of a result of
Zebulski and Zeevi (\cite{[ZZ]}). We will give here a simple proof for it  based on our
previous results.

\begin{theorem}
Let $g\in L^2({\Bbb R}^d)$ and $a = \frac{p}{q}$ be a rational number with $p<q$
and $p,q$ are co-prime. Define $g_j = Zg(x-\frac{p}{q}j,t)$ for $j\in\{0,1\cdots,q-1\}^d$.
If ${\mathcal G}(g,a,1)$ is a Gabor frame of $L^2({\Bbb R}^d)$, then there exists $A,B$ such that
\begin{equation}\label{eq5.3}
0<A\leq\max_{j\in\{0,\cdots,q-1\}^d}|Zg_j|\leq B<\infty \ \mbox{a.e. on $[0,1]^{2d}$.}
\end{equation}
If $a = \frac{1}{q}$, the converse also holds.
\end{theorem}

\begin{pf}
Note that $Z$ is a unitary mapping between $L^2({\Bbb R}^d)$ and $L^2([0,1]^{2d})$
and that ${\mathcal G}(g,a,1)$ is a Gabor frame on $L^2({\Bbb R}^d)$ if and only if
the image of the Gabor system under the Zak transform $Z[{\mathcal G}(g,a,1)]$ is a frame on $L^2([0,1]^{2d})$.
 Writing  $n = rq+j$ with $r\in{\Bbb Z}^d$ and $j\in\{0,\cdots,q-1\}^d$, we have
$$
\begin{aligned}
Z\left[e^{2\pi i \langle m,\cdot\rangle}g(\cdot-n\frac{p}{q})\right](x,t)
=& e^{2\pi i \langle m,x\rangle}\sum_{k\in{\Bbb Z}}g\left(x-k-pr-\frac{p}{q}j\right)
e^{2\pi i \langle k,t\rangle}\\
=&Zg(x-\frac{p}{q}j,t)e^{2\pi i \langle m,x\rangle}e^{2\pi i \langle rp,t\rangle}.
\end{aligned}
$$
From this, we see that
$$
\begin{aligned}
Z[{\mathcal G}(g,a,1)] =& \bigcup_{j\in\{0,1\cdots q-1\}^d}\bigcup_{m,r\in{\Bbb Z}^d}
\{g_j(x)e^{2\pi i \langle m,x\rangle}e^{2\pi i \langle rp,t\rangle}\}\\
=&\bigcup_{j\in\{0,1\cdots q-1\}^d}{\mathcal E}(g_j,{\Bbb Z}^d\times p{\Bbb Z}^d).
\end{aligned}
$$
i.e. $Z({\mathcal G}(g,a,1))$ is a system of  windowed exponentials on $[0,1]^{2d}$.
Therefore, if ${\mathcal G}(g,a,1)$ form a Gabor frame, then $Z[{\mathcal G}(g,a,1)]$
forms a frame of windowed exponentials. Moreover, ${\Bbb Z}^d\times p{\Bbb Z}^d$ has
positive upper Beurling density. By Theorem \ref{th3.1}, (\ref{eq5.3}) has to hold.

Conversely, if $a = 1/q$, then the exponential frequency set becomes ${\Bbb Z}^{2d}$
and the associated set of exponentials is an orthonormal basis for $L^2([0,1]^{2d})$.
According to the proof of Theorem \ref{th0.2}, (\ref{eq5.3}) implies that
$Z[{\mathcal G}(g,a,1)]$ forms a frame of windowed exponentials on $L^2([0,1]^{2d})$.
 Therefore, the original Gabor system forms a frame for $L^2({\Bbb R}^d)$.
This completes the proof.
\end{pf}

Zibulski and Zeevi \cite{[ZZ]} showed when $a=1/q$, ${\mathcal G}(g,a,1)$ is a
Gabor frame if and only if $\sum_{j=0}^{q-1}|g_j|^2$ is bounded above and
bounded away from $0$. Our result is consistent with their characterization since
$\ell^2$-norm and $\ell^{\infty}$-norm are equivalent on ${\Bbb R}^q$. For $a=p/q$,
our condition gives a simple necessary condition. For a necessary and sufficient
condition, we  refer the reader to Zibulski and Zeevi \cite{[ZZ]}, who expressed it in
terms of the  boundedness of the eigenvalues of an associated positive-definite matrix.

\end{document}